\documentclass[9pt, journal]{IEEEtran}

\usepackage{cite}
\usepackage{graphicx}

\usepackage{amsfonts}
\usepackage{amsmath} 
\usepackage{amssymb}

\usepackage{amsthm}	
\usepackage{enumerate}
\usepackage{mathrsfs}
\usepackage{subcaption}
\usepackage{epstopdf}
\usepackage{dashrule}

\newtheorem{lemm}{Lemma}
\newtheorem{definitio}{Definition}

\newtheorem{assumptio}{Assumption}

\newtheorem{propositio}{Proposition}

\newtheorem{remar}{Remark}

\def\nn{\nonumber}
\def\beq{\begin{equation}}
\def\eeq{\end{equation}}
\def\beqs{\begin{equation*}}
\def\eeqs{\end{equation*}}
\def\balign{\begin{align}}
\def\ealign{\end{align}}
\def\bea{\begin{eqnarray}}
\def\eea{\end{eqnarray}}
\def\ba{\begin{array}}
\def\ea{\end{array}}
\def\bcenter{\begin{center}}
\def\ecenter{\end{center}}
\def\bitem{\begin{itemize}}
\def\eitem{\end{itemize}}
\def\benum{\begin{enumerate}}
\def\eenum{\end{enumerate}}
\def\bdoc{\begin{document}}
\def\edoc{\end{document}}
\def\defeq{{\stackrel{\Delta}{=}}}
\def\nn{\nonumber}

\def\eg{{e.g., \/}}
\def\etal{{\it et al. \/}}
\def\viz{{\iet viz.,\ \/}}
\def\ie{{i.e.,\ \/}}
\def\etc{{\it etc.,\ \/}}
\def\cf{{cf.,\ \/}}

\usepackage[usenames,dvipsnames,svgnames,table]{xcolor}

\def\yellow{\textcolor{yellow}}
\def\green{\textcolor{green}}
\def\tcbg{\textcolor{bgrd}}
\def\magenta{\textcolor{magenta}}
\def\white{\textcolor{white}}
\def\gray{\textcolor{gray}}

\newcommand\Quote[1]{\lq\textsl{#1}\rq}
\newcommand\fr[2]{{\textstyle\frac{#1}{#2}}}

\newcommand{\Ambb}{\mathbb{A}}
\newcommand{\Bmbb}{\mathbb{B}}
\newcommand{\Cmbb}{\mathbb{C}}
\newcommand{\Dmbb}{\mathbb{D}}
\newcommand{\Embb}{\mathbb{E}}
\newcommand{\Fmbb}{\mathbb{F}}
\newcommand{\Gmbb}{\mathbb{G}}
\newcommand{\Hmbb}{\mathbb{H}}
\newcommand{\Imbb}{\mathbb{I}}
\newcommand{\Jmbb}{\mathbb{J}}
\newcommand{\Kmbb}{\mathbb{K}}
\newcommand{\Lmbb}{\mathbb{L}}
\newcommand{\Mmbb}{\mathbb{M}}
\newcommand{\Nmbb}{\mathbb{N}}
\newcommand{\Ombb}{\mathbb{O}}
\newcommand{\Pmbb}{\mathbb{P}}
\newcommand{\Qmbb}{\mathbb{Q}}
\newcommand{\Rmbb}{\mathbb{R}}
\newcommand{\Smbb}{\mathbb{S}}
\newcommand{\Tmbb}{\mathbb{T}}
\newcommand{\Umbb}{\mathbb{U}}
\newcommand{\Vmbb}{\mathbb{V}}
\newcommand{\Wmbb}{\mathbb{W}}
\newcommand{\Xmbb}{\mathbb{X}}
\newcommand{\Ymbb}{\mathbb{Y}}
\newcommand{\Zmbb}{\mathbb{Z}}

\newcommand{\Amsc}{\mathcal{A}}
\newcommand{\Bmsc}{\mathcal{B}}
\newcommand{\Cmsc}{\mathcal{C}}
\newcommand{\Dmsc}{\mathcal{D}}
\newcommand{\Emsc}{\mathcal{E}}
\newcommand{\Fmsc}{\mathcal{F}}
\newcommand{\Gmsc}{\mathcal{G}}
\newcommand{\Hmsc}{\mathcal{H}}
\newcommand{\Imsc}{\mathcal{I}}
\newcommand{\Jmsc}{\mathcal{J}}
\newcommand{\Kmsc}{\mathcal{K}}
\newcommand{\Lmsc}{\mathcal{L}}
\newcommand{\Mmsc}{\mathcal{M}}
\newcommand{\Nmsc}{\mathcal{N}}
\newcommand{\Omsc}{\mathcal{O}}
\newcommand{\Pmsc}{\mathcal{P}}
\newcommand{\Qmsc}{\mathcal{Q}}
\newcommand{\Rmsc}{\mathcal{R}}
\newcommand{\Smsc}{\mathcal{S}}
\newcommand{\Tmsc}{\mathcal{T}}
\newcommand{\Umsc}{\mathcal{U}}
\newcommand{\Vmsc}{\mathcal{V}}
\newcommand{\Wmsc}{\mathcal{W}}
\newcommand{\Xmsc}{\mathcal{X}}
\newcommand{\Ymsc}{\mathcal{Y}}
\newcommand{\Zmsc}{\mathcal{Z}}

\newcommand{\Ascr}{\mathscr{A}}
\newcommand{\Bscr}{\mathscr{B}}
\newcommand{\Cscr}{\mathscr{C}}
\newcommand{\Dscr}{\mathscr{D}}
\newcommand{\Escr}{\mathscr{E}}
\newcommand{\Fscr}{\mathscr{F}}
\newcommand{\Gscr}{\mathscr{G}}
\newcommand{\Hscr}{\mathscr{H}}
\newcommand{\Iscr}{\mathscr{I}}
\newcommand{\Jscr}{\mathscr{J}}
\newcommand{\Kscr}{\mathscr{K}}
\newcommand{\Lscr}{\mathscr{L}}
\newcommand{\Mscr}{\mathscr{M}}
\newcommand{\Nscr}{\mathscr{N}}
\newcommand{\Oscr}{\mathscr{O}}
\newcommand{\Pscr}{\mathscr{P}}
\newcommand{\Qscr}{\mathscr{Q}}
\newcommand{\Rscr}{\mathscr{R}}
\newcommand{\Sscr}{\mathscr{S}}
\newcommand{\Tscr}{\mathscr{T}}
\newcommand{\Uscr}{\mathscr{U}}
\newcommand{\Vscr}{\mathscr{V}}
\newcommand{\Wscr}{\mathscr{W}}
\newcommand{\Xscr}{\mathscr{X}}
\newcommand{\Yscr}{\mathscr{Y}}
\newcommand{\Zscr}{\mathscr{Z}}

\def\Atiny{\mbox{\tiny{A}}}
\def\Btiny{\mbox{\tiny{B}}}
\def\Ctiny{\mbox{\tiny{C}}}
\def\Dtiny{\mbox{\tiny{D}}}
\def\Etiny{\mbox{\tiny{E}}}
\def\Ctiny{\mbox{\tiny{C}}}
\def\Ftiny{\mbox{\tiny{F}}}
\def\Gtiny{\mbox{\tiny{G}}}
\def\Htiny{\mbox{\tiny{H}}}
\def\Itiny{\mbox{\tiny{I}}}
\def\Jtiny{\mbox{\tiny{J}}}
\def\Ktiny{\mbox{\tiny{K}}}
\def\Ltiny{\mbox{\tiny{L}}}
\def\Mtiny{\mbox{\tiny{M}}}
\def\Ntiny{\mbox{\tiny{N}}}
\def\Otiny{\mbox{\tiny{O}}}
\def\Ptiny{\mbox{\tiny{P}}}
\def\Qtiny{\mbox{\tiny{Q}}}
\def\Rtiny{\mbox{\tiny{R}}}
\def\Stiny{\mbox{\tiny{S}}}
\def\Ttiny{\mbox{\tiny{T}}}
\def\Utiny{\mbox{\tiny{U}}}
\def\Vtiny{\mbox{\tiny{V}}}
\def\Wtiny{\mbox{\tiny{W}}}
\def\Xtiny{\mbox{\tiny{X}}}
\def\Ytiny{\mbox{\tiny{Y}}}
\def\Ztiny{\mbox{\tiny{Z}}}

\def\alphabf{\hbox{\boldmath$\alpha$\unboldmath}}
\def\betabf{\hbox{\boldmath$\beta$\unboldmath}}
\def\gammabf{\hbox{\boldmath$\gamma$\unboldmath}}
\def\deltabf{\hbox{\boldmath$\delta$\unboldmath}}
\def\epsilonbf{\hbox{\boldmath$\epsilon$\unboldmath}}
\def\zetabf{\hbox{\boldmath$\zeta$\unboldmath}}
\def\etabf{\hbox{\boldmath$\eta$\unboldmath}}
\def\iotabf{\hbox{\boldmath$\iota$\unboldmath}}
\def\kappabf{\hbox{\boldmath$\kappa$\unboldmath}}
\def\lambdabf{\hbox{\boldmath$\lambda$\unboldmath}}
\def\mubf{\hbox{\boldmath$\mu$\unboldmath}}
\def\nubf{\hbox{\boldmath$\nu$\unboldmath}}
\def\xibf{\hbox{\boldmath$\xi$\unboldmath}}
\def\pibf{\hbox{\boldmath$\pi$\unboldmath}}
\def\rhobf{\hbox{\boldmath$\rho$\unboldmath}}
\def\sigmabf{\hbox{\boldmath$\sigma$\unboldmath}}
\def\taubf{\hbox{\boldmath$\tau$\unboldmath}}
\def\upsilonbf{\hbox{\boldmath$\upsilon$\unboldmath}}
\def\phibf{\hbox{\boldmath$\phi$\unboldmath}}
\def\chibf{\hbox{\boldmath$\chi$\unboldmath}}
\def\psibf{\hbox{\boldmath$\psi$\unboldmath}}
\def\omegabf{\hbox{\boldmath$\omega$\unboldmath}}
\def\Sigmabf{\hbox{$\bf \Sigma$}}
\def\Upsilonbf{\hbox{$\bf \Upsilon$}}
\def\Omegabf{\hbox{$\bf \Omega$}}
\def\Deltabf{\hbox{$\bf \Delta$}}
\def\Gammabf{\hbox{$\bf \Gamma$}}
\def\Thetabf{\hbox{$\bf \Theta$}}
\def\Lambdabf{\mbox{$ \bf \Lambda $}}
\def\Xibf{\hbox{\bf$\Xi$}}
\def\Pibf{{\bf \Pi}}
\def\thetabf{{\mbox{\boldmath$\theta$\unboldmath}}}
\def\Upsilonbf{\hbox{\boldmath$\Upsilon$\unboldmath}}
\newcommand{\Phibf}{\mbox{${\bf \Phi}$}}
\newcommand{\Psibf}{\mbox{${\bf \Psi}$}}

\def\abf{{\bf a}}
\def\bbf{{\bf b}}
\def\cbf{{\bf c}}
\def\dbf{{\bf d}}
\def\ebf{{\bf e}}
\def\fbf{{\bf f}}
\def\gbf{{\bf g}}
\def\hbf{{\bf h}}
\def\ibf{{\bf i}}
\def\jbf{{\bf j}}
\def\kbf{{\bf k}}
\def\lbf{{\bf l}}
\def\mbf{{\bf m}}
\def\nbf{{\bf n}}
\def\obf{{\bf o}}
\def\pbf{{\bf p}}
\def\qbf{{\bf q}}
\def\rbf{{\bf r}}
\def\sbf{{\bf s}}
\def\tbf{{\bf t}}
\def\ubf{{\bf u}}
\def\vbf{{\bf v}}
\def\wbf{{\bf w}}
\def\xbf{{\bf x}}
\def\ybf{{\bf y}}
\def\zbf{{\bf z}}
\def\rbf{{\bf r}}
\def\xbf{{\bf x}}
\def\ybf{{\bf y}}
\def\Abf{{\bf A}}
\def\Bbf{{\bf B}}
\def\Cbf{{\bf C}}
\def\Dbf{{\bf D}}
\def\Ebf{{\bf E}}
\def\Fbf{{\bf F}}
\def\Gbf{{\bf G}}
\def\Hbf{{\bf H}}
\def\Ibf{{\bf I}}
\def\Jbf{{\bf J}}
\def\Kbf{{\bf K}}
\def\Lbf{{\bf L}}
\def\Mbf{{\bf M}}
\def\Nbf{{\bf N}}
\def\Obf{{\bf O}}
\def\Pbf{{\bf P}}
\def\Qbf{{\bf Q}}
\def\Rbf{{\bf R}}
\def\Sbf{{\bf S}}
\def\Tbf{{\bf T}}
\def\Ubf{{\bf U}}
\def\Vbf{{\bf V}}
\def\Wbf{{\bf W}}
\def\Xbf{{\bf X}}
\def\Ybf{{\bf Y}}
\def\Zbf{{\bf Z}}
\def\Ac{{\cal A}}
\def\Bc{{\cal B}}
\def\Cc{{\cal C}}
\def\Dc{{\cal D}}
\def\Ec{{\cal E}}
\def\Fc{{\cal F}}
\def\Gc{{\cal G}}
\def\Hc{{\cal H}}
\def\Ic{{\cal I}}
\def\Jc{{\cal J}}
\def\Kc{{\cal K}}
\def\Lc{{\cal L}}
\def\Mc{{\cal M}}
\def\Nc{{\cal N}}
\def\Oc{{\cal O}}
\def\Pc{{\cal P}}
\def\Qc{{\cal Q}}
\def\Rc{{\cal R}}
\def\Sc{{\cal S}}
\def\Tc{{\cal T}}
\def\Uc{{\cal U}}
\def\Vc{{\cal V}}
\def\Wc{{\cal W}}
\def\Xc{{\cal X}}
\def\Yc{{\cal Y}}
\def\Zc{{\cal Z}}

\def\0bf{\mathbf{0}}
\def\1bf{\mathbf{1}}

\def\lambdabf{\boldsymbol{\lambda}}
\def\mubf{\boldsymbol{\mu}}
\def\gammabf{\boldsymbol{\gamma}}
\def\xibf{\boldsymbol{\xi}}
\def\sigmabf{\boldsymbol{\sigma}}
\def\zetabf{\boldsymbol{\zeta}}

\def\Lambdabf{\boldsymbol{\Lambda}}
\def\Mubf{\boldsymbol{\Mu}}
\def\Gammabf{\boldsymbol{\Gamma}}
\def\Xibf{\boldsymbol{\Xi}}
\def\Sigmabf{\boldsymbol{\Sigma}}
\def\Zetabf{\boldsymbol{\Zeta}}

\newcommand\independent{\protect\mathpalette{\protect\independenT}{\perp}}
\def\independenT#1#2{\mathrel{\rlap{$#1#2$}\mkern2mu{#1#2}}}

\newcommand{\bsf}[1]{\textsf{\textbf{#1}}}
\newcommand{\lbsf}[1]{\textsf{\large  \textbf{#1}}}
\newcommand{\Lbsf}[1]{\textsf{\Large  \textbf{#1}}}
\newcommand{\hbsf}[1]{\textsf{\huge  \textbf{#1}}}

\newcommand{\myminipage}[3]{\begin{minipage}[#1]{#2}{#3} \end{minipage}}
\newcommand{\sbs}[4]{\myminipage{c}{#1}{#3} \hfill
\myminipage{c}{#2}{#4}}

\newcommand{\myfig}[2]{\centerline{\psfig{figure=#1,width=#2,silent=}}}
\newcommand{\myfigh}[2]{\centerline{\psfig{figure=#1,height=#2,silent=}}}
\newcommand{\myfigwh}[3]{\centerline{\psfig{figure=#1,width=#2,height=#3,silent=}}}

\newcommand{\beqa}{\begin{eqnarray}}
\newcommand{\eeqa}{\end{eqnarray}}
\newcommand{\beqan}{\begin{eqnarray*}}
\newcommand{\eeqan}{\end{eqnarray*}}
\newcommand{\dst}[1]{\displaystyle{ #1 }}
\newcounter{pf}

\newcommand{\smax}[1] { \bar \sigma \left( #1 \right) }
\newcommand{\Rn}{{\mathbb R}^n}
\newcommand{\R}{{\mathbb R}}
\newcommand{\C}{{\mathbb C}}
\newcommand{\Rm}{\mathbb{R}^m}
\newcommand{\Rmn}{\mathbb{R}^{m \times n}}
\newcommand{\Rpq}{\mathbb{R}^{p \times q}}
\newcommand{\Cn}{\mathbb{C}^n}
\newcommand{\Cm}{\mathbb{C}^m}
\newcommand{\Cnn}{\mathbb{C}^{n \times n}}
\newcommand{\Cmn}{\mathbb{C}^{m \times n}}
\newcommand{\ip}[1]{\left\langle #1 \right\rangle}
\newcommand{\rank}{\mbox{rank}}
\newcommand{\Span}{\mbox{\rm Span }}
\newcommand{\Trace}{\mbox{\rm Tr }}
\newcommand{\trace}[1]{\text{Tr}\left(#1\right)}
\newcommand{\Spec}{\mbox{\rm Spec }}
\newcommand{\vectornorm}[1]{\left\|#1\right\|}

\newcommand{\pd}[2]{\frac{\partial #1}{\partial #2}}
\newcommand{\ppd}[3]{\frac{\partial^2 #1}{\partial #2 \partial #3}}

\newcommand{\thtilde}{\tilde{\theta}}
\newcommand{\thnom}{\theta^\circ}
\newcommand{\thopt}{\theta^{\mbox{\small opt}}}
\newcommand{\thhat}{{\hat{\theta}}}
\newcommand{\Tho}{\Theta^\circ}
\newcommand{\tho}{\theta^\circ}
\newcommand{\np}{{n_p}}

\newcommand{\ii}{{[i]}}
\newcommand{\II}{{[i+1]}}
\newcommand{\iii}{{[ii]}}
\newcommand{\jj}{{[j]}}
\newcommand{\kk}{{[k]}}
\newcommand{\thi}{{\theta^\ii}}
\newcommand{\thI}{{\theta^\II}}
\newcommand{\di}{{d^\ii}}
\newcommand{\gi}{{g^\ii}}
\newcommand{\Hi}{{\HH^\ii}}
\newcommand{\thK}{\theta^{(k+1)}}
\newcommand{\gk}{{g^{(k)}}}
\newcommand{\Hk}{{{\cal H}^{(k)}}}

\newcommand{\bfdelta}{{\bf \Delta}}

\newcommand{\Exp}[1]{\exp \left\{ #1 \right\}} 
\newcommand{\gaussian}[1]{\mathbb{N} \left( #1 \right)}
\newcommand{\uniform}[1]{\mathbb{U} \left[ #1 \right]}
\newcommand{\exponential}[1]{\mathbb{E} \left[ #1 \right]}
\newcommand{\EXP}[1]{\EEXP \left[ #1 \right]} 
\newcommand{\EEXP}{\mbox{\bsf{E}}} 
\newcommand{\Prob}[1]{\mbox{{\sf Pr}} \left(#1 \right)}
\newcommand{\convas}{\stackrel{as}{\longrightarrow}}
\newcommand{\convinp}{\stackrel{p}{\longrightarrow}}
\newcommand{\convind}{\stackrel{d}{\longrightarrow}}
\newcommand{\convqm}{\stackrel{qm}{\longrightarrow}}
\newcommand{\sss}[1]{{_{#1}}}
\newcommand{\density}[2]{p_{_{_{#1}}}\!\!\left(#2 \right)} 
\newcommand{\distro}[2]{P_{_{_{#1}}}\!\!\left(#2 \right)} 
\newcommand{\rxx}[1]{R_{_{#1}}\!} 
\newcommand{\sxx}[1]{S_{_{#1}}} 
\newcommand{\cov}[1]{\Lambda_{_{#1}}} 
\newcommand{\mean}[1]{m_{_{#1}}} 
\newcommand{\LS}[1]{\hat{#1}_{_{LS}}} 
\newcommand{\MV}[1]{\hat{#1}_{_{MV}}} 
\newcommand{\LMV}[1]{\hat{#1}_{_{LMV}}} 
\newcommand{\ML}[1]{\hat{#1}_{_{ML}}} 

\renewcommand{\arraystretch}{0.9}
\newcommand{\bmat}[1]{ \begin{bmatrix} #1 \end{bmatrix}}
\newcommand{\mat}[1]{ \left[ \begin{array}{cccccccc} #1 \end{array}
\right] }
\newcommand{\smallmat}[1]{\small{\mat{#1}}}
\newcommand{\sysblk}[4]{\begin{array}{c|cccc}#1&#2\\ \hline#3&#4
\end{array}}
\newcommand{\sysmat}[4]{\left[\sysblk{#1}{#2}{#3}{#4}\right]}
\newcommand{\SGeq}{\succ}
\newcommand{\SLeq}{\prec}
\newcommand{\Geq}{\succeq}
\newcommand{\Leq}{\preceq}

\newcommand{\Aset}{\mathbb{A}}
\newcommand{\Bset}{\mathbb{B}}
\newcommand{\Fset}{\mathbb{F}}
\newcommand{\Gset}{\mathbb{G}}
\newcommand{\Kset}{\mathbb{K}}
\newcommand{\Mset}{\mathbb{M}}
\newcommand{\Sset}{\mathbb{S}}
\newcommand{\Tset}{\mathbb{T}}
\newcommand{\Uset}{\mathbb{U}}
\newcommand{\Vset}{\mathbb{V}}
\newcommand{\Wset}{\mathbb{W}}

\newcommand{\Ical}{{\cal I}}
\newcommand{\Acal}{{\cal A}}
\newcommand{\Bcal}{{\cal B}}
\newcommand{\Ccal}{{\mathcal{C}}}
\newcommand{\Dcal}{{\cal D}}
\newcommand{\Ecal}{{\mathcal{E}}}
\newcommand{\Fcal}{{\cal F}}
\newcommand{\Gcal}{{\mathcal{G}}}
\newcommand{\Hcal}{{\cal H}}
\newcommand{\Kcal}{{\mathcal{K}}}
\newcommand{\Lcal}{{\cal L}}
\newcommand{\Mcal}{{\cal M}}
\newcommand{\Ncal}{{\mathcal{N}}}
\newcommand{\Pcal}{{\cal P}}
\newcommand{\Qcal}{{\mathcal{Q}}}
\newcommand{\Rcal}{{\cal R}}
\newcommand{\Scal}{{\mathcal{S}}}
\newcommand{\Tcal}{{\mathcal{T}}}
\newcommand{\Wcal}{{\mathcal{W}}}
\newcommand{\Ucal}{{\cal U}}
\newcommand{\Vcal}{{\mathcal{V}}}
\newcommand{\Xcal}{{\cal X}}
\newcommand{\Ycal}{{\mathcal{Y}}}
\newcommand{\Zcal}{{\mathcal{Z}}}

\newcommand{\EE}{{\bf E}}
\newcommand{\FF}{{\bf F}}
\newcommand{\GG}{{\bf G}}
\newcommand{\HH}{{\bf H}}
\newcommand{\LL}{{\bf L}}
\newcommand{\NN}{{\bf N}}
\newcommand{\MM}{{\bf M}}
\newcommand{\PP}{{\bf P}}
\newcommand{\QQ}{{\bf Q}}
\newcommand{\RR}{{\bf R}}
\renewcommand{\SS}{{\bf S}}
\newcommand{\TT}{{\bf T}}
\newcommand{\VV}{{\bf V}}
\newcommand{\WW}{{\bf W}}

\newcommand{\thk}{\theta^{(k)}}
\newcommand{\thb}{\theta^{\rm opt}}
\newcommand{\alb}{\alpha^{\rm opt}}
\newcommand{\dk}{d^{(k)}}
\newcommand{\Hinf}{{\cal H}_\infty}
\newcommand{\Htwo}{{\cal H}_2}

\renewcommand{\arraystretch}{1.1}

\newcommand{\red}[1]{{\color{red} #1}}
\newcommand{\blue}[1]{{\color{Blue} #1}}
\newcommand{\black}[1]{{\color{Black} #1}}

\newcounter{l1}
\newcounter{l2}
\newcounter{l3}
\newcommand{\bdotlist}{\begin{list}{$\bullet$}{}}
\newcommand{\bboxlist}{\begin{list}{$\Box$}{}}
\newcommand{\bbboxlist}{\begin{list}{\raisebox{.005in}{{\tiny
$\blacksquare$ \ \ }}}{}}
\newcommand{\bdashlist}{\begin{list}{$-$}{} }
\newcommand{\blist}{\begin{list}{}{} }
\newcommand{\barablist}{\begin{list}{\arabic{l1}}{\usecounter{l1}}}
\newcommand{\balphlist}{\begin{list}{(\alph{l2})}{\usecounter{l2}}}
\newcommand{\bAlphlist}{\begin{list}{\Alph{l2}.}{\usecounter{l2}}}
\newcommand{\bdiamlist}{\begin{list}{$\diamond$}{}}
\newcommand{\bromalist}{\begin{list}{(\roman{l3})}{\usecounter{l3}}}

\newcommand{\prf}[1]{  \noindent {\em Proof.} \;  #1 \hfill $\blacksquare$}
\usepackage{tikz}
\usetikzlibrary{arrows,snakes,backgrounds}
\usepackage{pgfplots}
\pgfplotsset{compat=newest}
\usetikzlibrary{plotmarks}
\usepackage{grffile}

\newcommand{\argmin}{\mathop{\rm argmin}}
\newcommand{\argmax}{\mathop{\rm argmax}}
\newcommand{\diag}{\mathop{\mathrm{diag}}}
\newcommand{\tr}{\mathop{\rm Tr}}
\newcommand{\conv}{\mathop{\rm conv}}
\newcommand{\var}{\mathop{\rm Var}}
\renewcommand{\b}[1]{\ensuremath{\boldsymbol{\mathrm{#1}}}}
\newcommand{\Asur}{{\widetilde{A}}}
\newcommand{\Bsur}{{\widetilde{B}}}
\newcommand{\Gsur}{{\widetilde{L}}}
\newcommand{\Hsur}{{\widetilde{H}}}

\newcommand{\E}[1]{\b{\mu}_{{#1}}}
\newcommand{\Var}[1]{{\Sigma_{#1}}}

\newcommand{\bone}{\mathbf{1}}

\def\fic{v}
\def\contract{\Vcal}
\def\fictransform{Z}
\def\surrpol{\phi}
\def\surrpolset{\Phi}
\def\Xw{{P^w}}
\def\Xxi{{P^\xi}}
\def \couple{\Ccal}
\def \nodes{V}
\def \edges{E}
\def \gr{G}
\def \gnest{\gr_N}
\def \gcouple{\gr_C}
\def \enest{\edges_N}
\def \ecouple{\edges_C}
\def\support{\Wcal}
\def \pinfo{\zeta}
\def\proj{\Pi_\couple}
\def \gy{\gr}
\def \ey{\edges}
\def \gyset{\Gcal}
\def \lmi{\sf LMI}
\def \wtoX{L}

\def \add [#1]{\blue{#1}}
\def \remove [#1]{\red{#1}}
\def \replace [#1]#2{\red{#1} \blue{#2}}

\newcommand{\wl}[1]{\blue{(\textbf{Sam says:} #1)}}
\newcommand{\eb}[1]{\red{(\textbf{Eilyan says:} #1)}}

\newcommand{\rev}[1]{#1}
\newcommand{\conf}[1]{#1}
\newcommand{\rauto}[1]{#1}
\newcommand{\rtwo}[1]{#1}
\newcommand{\cut}[1]{\red{#1}}
\newcommand{\sam}[1]{#1}

\graphicspath{{figures/}}

\begin{document}

\title{Design of Robust Decentralized Controllers via Assume-Guarantee Contracts
}

\author{\vspace{.07in} Weixuan Lin   \qquad Eilyan Bitar 
\thanks{Supported in part by the National Science Foundation under grants ECCS-135162 and IIP-1632124, and the Holland Sustainability Project Trust. This paper builds on the authors'
preliminary results published as part of the 2020 IEEE American Control Conference \cite{lin2020decentralized}. This paper differs from the conference version in terms of formal mathematical proofs omitted from the conference version, and expanded technical discussions and remarks throughout the paper. }
\thanks{Weixuan Lin ({\tt\small wl476@cornell.edu}) and Eilyan Bitar ({\tt\small eyb5@cornell.edu}) are with the School of Electrical and Computer Engineering, Cornell University, Ithaca, NY, 14853, USA. }
\vspace{-.25in}
}

\maketitle

\begin{abstract}
We consider the decentralized control of a discrete-time time-varying linear system subject to additive disturbances and polyhedral constraints on the state and input trajectories. The underlying system is composed of a finite collection of  dynamically coupled subsystems, where each subsystem is assumed to have a dedicated local controller. The decentralization of information is expressed according to sparsity constraints on the state measurements that each local controller has access to. We investigate the design of decentralized controllers that are affinely parameterized in their measurement history. For problems with partially nested information structures, the optimization over such restricted policy spaces is known to be convex. Convexity is not, however, guaranteed under more general (nonclassical) information structures, where  the information available to one local controller can be affected by control actions that it cannot  reconstruct. To address the nonconvexity that arises in such problems, we propose an  approach to decentralized control design where such information-coupling  states  are effectively treated as disturbances whose trajectories are constrained to take values in ellipsoidal ``contract'' sets whose location, scale, and orientation are  jointly optimized with the underlying affine decentralized control policy. We establish a  structural condition on the space of allowable contracts that facilitates the joint optimization over the control policy and the contract set via semidefinite programming.
\end{abstract}



\section{Introduction}

We investigate the design of affine decentralized control policies for 
stochastic discrete-time, linear systems that evolve over a finite horizon, and are subject to polyhedral constraints on the state and input trajectories. The computational tractability of such problems depends in part on their information structures  \cite{sandell1974solution, Tsitsiklis1985}. 
A decentralized control problem is said to have a \emph{nonclassical} information structure if the information available to one controller can be affected by the control actions of another that it cannot access or reconstruct. 
Under such information structures, the calculation of optimal decentralized control policies is known to be computationally intractable, because of the so called signalling incentive for controllers to communicate with each other via the actions they undertake \cite{Witsenhausen1968, sandell1974solution, Tsitsiklis1985}.
To complicate matters further, there may be hard constraints  coupling the local actions and states of different controllers that must be jointly enforced without explicit communication.
In this paper, we  address these challenges by relaxing the requirement that decentralized controllers be optimal with respect to the broad family of all causal policies, and instead search for suboptimal decentralized controllers that can be efficiently computed via convex optimization.

\emph{Related Literature: \ }
\rev{Many decentralized control methods based on  model predictive control (MPC) tools have been proposed \cite{venkat2005stability, Camponogara2002, Jia2001, Bemporad2010, Keviczky2006, dunbar2007distributed, Dunbar2006, lucia2015contract, farina2012distributed, richards2004decentralized, richards2007robust, riverso2012tube, trodden2006robust, trodden2010distributed, hernandez2017distributed, trodden2014cooperative}.}
These approaches typically rely on a decomposition of the decentralized control problem into a collection of decoupled  local control problems, where each subproblem is solved using traditional MPC methods.
A particular class of methods---commonly referred to as ``tube-based'' decentralized MPC---treat the coupling states and inputs affecting each subsystem as independent and bounded disturbances that are assumed to take values in the given state and input constraint sets  \cite{dunbar2007distributed, Dunbar2006, lucia2015contract, farina2012distributed, richards2004decentralized, richards2007robust, riverso2012tube, trodden2006robust, trodden2010distributed, hernandez2017distributed, trodden2014cooperative}. 
Decentralized controllers  computed in this manner  may result in overly conservative behaviors for a
number of reasons.
\emph{First}, the treatment of the coupling states and inputs as independent disturbances ignores the potential dynamical coupling between these variables.
\emph{Second}, the over approximation of the coupling state and input trajectory sets by their corresponding state and input constraint sets will likely be  loose for many problem instances. More importantly, the over approximation of the coupling state and input trajectory sets in this manner ignores the fact that these sets depend on the control policy being used to regulate the system, and, therefore, neglects the possibility of co-optimizing their specification with the control policy.

\emph{Contribution: \ } 
We provide a computationally tractable method to calculate control policies that are guaranteed to be feasible for constrained decentralized control problems with nonclassical information structures.
Loosely speaking, the proposed approach  eliminates the informational coupling between subsystems by treating the information-coupling states as  disturbances whose trajectories are ``assumed" to take values in  a  ``contract'' set that we design. To ensure the satisfaction of this assumption, we impose a contractual constraint on the control policy that ``guarantees'' that the information-coupling states that it induces indeed  belong to the contract set.
Naturally, this approach yields an inner approximation of the original decentralized control design problem, where the conservatism of the resulting approximation depends on the specification of the contract set.
To limit the suboptimality of such approximations, we formulate a semi-infinite program to co-optimize the decentralized control policy with the location, scale, and orientation of an ellipsoidal contract set. \rtwo{The ability to adjust both the scale and orientation of the contract set enables the construction of contract sets that accurately capture the spatial and inter-temporal coupling between different subsystem states.} We establish a condition on the set of allowable contracts that facilitates the joint optimization of the control policy and the contract set via  semidefinite programming.

Assume-guarantee contracts  have  been utilized in a  variety of control applications \cite{kim2015compositional, chen2019compositional,nuzzo2013contract,sangiovanni2012taming}.
Closer to the setting considered in this paper, there are related papers in the literature that investigate an approach to decentralized control design via the co-optimization of the control policy and the contract sets that govern the interactions between the different subsystems in the network \cite{trodden2017distributed, darivianakis2018decentralized}. \rtwo{The techniques developed in these papers have two primary limitations.
First, they treat all states that result in a physical coupling between different subsystems as  disturbances, while our approach only treats states that induce an informational coupling between subsystems as disturbances. 
Second, these techniques only permit the scaling and translation of a base contract set when co-optimizing it with the control policy. To the best of our knowledge, the method proposed in this paper provides the first computationally tractable approach to co-optimizing the control policy with the location, scale, and \emph{orientation} of the contract set,  expanding substantially the family of contracts that can be efficiently optimized over. 
We provide an example that illustrates  the advantages of the proposed method.}

\emph{Notation: \ }   Let $\RR$ and $\RR_+$ denote the sets of real and non-negative real numbers, respectively.
Given a collection of vectors $x_1, .., x_N$ and an index set $J \subseteq \{1, .., N\}$, 
we let $x_J$ denote the concatenation of the vectors $x_j$  for $j \in J$ in ascending order of their indices. Given a sequence $\{x(t)\}$ and time indices $s \leq t$, we denote the sequence of elements from time $s$ to $t$ by  $x^{s:t} := (x(s), x(s+1), .., x(t) )$.
 Given a block matrix $A$, we let $[A]_{ij}$ denote its $(i,j)$-th block. We denote the trace of a square matrix $A$ by $\trace{A}$.
We denote the Minkowski sum of two sets $\Scal, \Tcal \subseteq \RR^n$  by $\Scal \oplus \Tcal := \{x + y \,  | \, x \in \Scal, \  y \in \Tcal \}$.

\section{Problem Formulation} \label{sec:formulation}

\subsection{System Model}

Consider a discrete-time, linear time-varying system consisting of $N$ coupled subsystems whose dynamics are described by
\begin{align} \label{eq:x_i}
x_i(t+1) = \sum_{j=1}^N  \left( A_{ij}(t) x_j(t)   + B_{ij} (t) u_j (t) \right) +  w_i (t),
\end{align}
for $i = 1, \dots , N$. 
We denote the \emph{local state}, \emph{local input}, and \emph{local disturbance} associated with  each subsystem $i$ at time $t$ by $x_i(t) \in \RR^{n_x^i}$, $u_i(t) \in \RR^{n_u^i}$, and $w_i (t) \in \RR^{n_x^i}$, respectively.
The system is assumed to evolve  over a finite time horizon $T$, and  the initial condition is assumed to be a random vector with known probability distribution.
Eq. \eqref{eq:x_i} can be presented more compactly as
\begin{align} \label{eq:full_sys}
x(t+1) = A(t) x(t) + B(t) u (t) + w(t).
\end{align}
Here, we denote by $x(t) := (x_1(t), .., x_N(t)) \in \RR^{n_x}$, $u(t) := (u_1(t), .., u_N(t)) \in \RR^{n_u}$, and $w (t) := (w_1(t), .., w_N(t)) \in \RR^{n_x}$ the \emph{full system state, input,} and \emph{disturbance} at time $t$. The dimensions of the system state and input are given by $n_x := \sum_{i=1}^N n_x^i$ and $n_u := \sum_{i=1}^N n_u^i$, respectively. 
The input and disturbance trajectories are related to the state trajectory according to
\begin{align}
x = B u + \wtoX w, \label{eq:trajectory}
\end{align}
where $x$, $u$, and $w$ denote the system state, input, and disturbance trajectories, respectively. They are defined according to
\begin{align}
x & := (x(0), \dots, x(T))  \in \RR^{N_x},  \label{eq:x_traj_def}\\ 
u & := (u(0), \dots, u(T-1))  \in \RR^{N_u}, \label{eq:u_traj_def}\\
w & : = (w(-1), w(0), \dots, w(T-1))  \in \RR^{N_x}, \label{eq:w_traj_def}
\end{align}
where the corresponding dimensions are given by 
$N_x := n_x(T+1)$ and $N_u := n_u T$.
To simplify the specification of disturbance-feedback affine control policies in the sequel, we adopt the notational convention $w (-1 ) := x(0)$. \rtwo{The lower block triangular 
matrices $B$ and $L$  are straightforward to  construct from the problem data in  \eqref{eq:full_sys}.}

\subsection{Disturbance Model}

We assume that disturbance trajectory $w$ is a zero-mean random vector. We denote its second moment matrix  by $M := \EE [w w^\top].$
We let  $\support \subset \RR^{N_x}$ denote the support of $w$, which we assume is a convex, compact set  with a non-empty interior. This assumption   ensures that the matrix $M$ is positive definite and  finite valued.

\subsection{System Constraints}

We consider a general family of polyhedral constraints on the state and input trajectories of the form 
\begin{align}
F_x x + F_u u + F_w w \leq g \qquad  \forall w \in \support, \label{eq:robust_con}
\end{align} 
where $F_x \in \RR^{m \times N_x}$, $F_u \in \RR^{m \times N_u}$, $F_w \in \RR^{m \times N_x}$, $g \in \RR^{m}$ are assumed to be given. Note that, in general, such constraints may couple states and inputs across different subsystems and time periods.

\subsection{Information Structure}

We encode the pattern according to which information is shared between subsystems with a directed graph $\gr_I = (\nodes,\edges_I)$, which we refer to as the \emph{information graph} of the system. 
The vertex set $\nodes = \{1, \dots , N \}$ indexes the subsystems being controlled, and the edge set $\edges_I \subseteq \nodes \times \nodes$ includes the directed edge
$(i,j)$ if and only if  subsystem $j$ has access to subsystem $i$'s local state measurement at each time $t$. 
We let $\nodes_I^-(i)$ denote the in-neighborhood of each subsystem $i \in \nodes$ in the information graph $\gr_I$.

Each subsystem is assumed to have access to the entire history of its local information up to and including time $t$.
Formally, the \emph{local information} available to each subsystem $i$ at time $t$ is defined as
\begin{align} \label{eq:pol1}
z_i(t) : = \{ x_j^{0:t} \ | \ (j, i) \in \edges_I \}.
\end{align}
The local control input to each subsystem $i$ is restricted to be a causal function of its local information. That is, the local input to subsystem $i$ at time $t$ is of the form
\begin{align} \label{eq:pol2}
u_i(t)  = \gamma_i(z_i(t), t),
\end{align}
where $\gamma_i(\cdot, t)$ is a measurable function of the local information $z_i (t)$.
We  define the \emph{local control policy} for subsystem $i$ as $\gamma_i : = (\gamma_i(\cdot,0), \dots, \gamma_i(\cdot,T-1))$. We refer to the  collection of local control policies  $\gamma : = (\gamma_1, \dots, \gamma_N)$ as the \emph{decentralized control policy}, which relates the state trajectory $x$ to the input trajectory $u$ according to $u = \gamma (x)$.
Finally, we let $\Gamma$ denote the set of all decentralized control policies respecting the information constraints encoded in Eq. \eqref{eq:pol2}.

\subsection{Decentralized Control Design}
We consider the following family of constrained decentralized control design problems:
\begin{align}
\begin{alignedat}{8}
&\text{minimize} \quad && \EE \left[ x^\top R_x x + u ^\top R_u u \right] \\
& \text{subject to} \quad && \gamma \in \Gamma  \\
&&& \hspace{-.105in}\left. \begin{array}{l}
u = \gamma (x) \\
x = B u + \wtoX w  \\
F_x x + F_u u + F_w w \leq g 
\end{array}
\hspace{-0.04in}\right\} \forall w \in \support.
\end{alignedat} \label{opt:decent}
\end{align}
Here, the cost matrices $R_x \in \RR^{N_x \times N_x}$ and $R_u \in \RR^{N_u \times N_u}$ are  assumed to be symmetric and positive semidefinite.
The tractability of the decentralized control design problem \eqref{opt:decent} depends  on its information structure.
In particular, if the information structure is \emph{partially nested}, then problem \eqref{opt:decent} can be equivalently reformulated (via the Youla parameterization) as a convex program in the space of disturbance feedback policies \cite{Ho1972}.
If the information structure is \emph{nonclassical} (i.e., not partially nested), then problem \eqref{opt:decent} is known to be computationally intractable, in general \cite{Tsitsiklis1985, Sandell1978, Mahajan2012_survey}.

\section{Information Decomposition} \label{sec:preliminaries}
The primary difficulty in solving decentralized control design problems stems from  the informational coupling that emerges when a subsystem's local information is affected by prior control actions that it cannot access or reconstruct. With the aim of isolating the effects of these actions on the  information available to each subsystem, we propose an   information  decomposition that partitions the local information available to each subsystem into a partially nested  subset (i.e., an information subset that is unaffected by control actions previously applied to the system) and its complement.
This  decomposition enables an equivalent reformulation of the decentralized control design problem where the control policy is expressed as an explicit function of the system disturbance and the so called \emph{information-coupling} states.
This reformulation will serve as the basis for the contract-based approach to decentralized control design  in Sec. \ref{sec:feas_control}.

\subsection{Decomposition of Local Information}
We decompose the local information available to each subsystem according to a partition of  its in-neighbors in the  information graph $\gr_I$.  Specifically, for each subsystem $i \in \nodes$, we let $\Ncal(i) \subseteq \nodes_I^{-}(i)$ denote the set of in-neighboring subsystems such that the information conveyed by their local state measurements is unaffected by the prior control actions of any subsystem. 
This requirement is satisfied if the local information of subsystem $i$  permits the reconstruction of all states and control actions directly affecting the local states of all subsystems belonging to $\Ncal(i)$. 
We denote the  complement of this set by $\couple(i) : = \nodes_I^{-}(i) \setminus \Ncal(i)$ for each subsystem $i \in \nodes$.

With the goal of providing an explicit characterization of this partition, we first provide a characterization of the physical coupling between different subsystems as reflected by the block sparsity patterns of the system matrices $A$ and $B$. We describe this coupling in terms of a pair of directed graphs, $\gr_A := (\nodes, \edges_A)$ and $\gr_B : = (\nodes, \edges_B)$, whose edge sets  are defined according to
\begin{align*}
\edges_A := \{ (j,i) \in \nodes \times \nodes \ | \  \exists t = 0, \dots, T-1 \text{ s.t. } A_{ij} (t)\neq 0\},\\
\edges_B := \{ (j,i) \in \nodes \times \nodes \ | \  \exists t = 0, \dots, T-1 \text{ s.t. } B_{ij} (t)\neq 0\}.
\end{align*}
We let  $\nodes_A^- (i)$ and $\nodes_B^- (i)$ denote the in-neighborhoods associated with each node $i \in \nodes$ in  $\gr_A$ and $\gr_B$, respectively.

Building on this graphical representation of the physical coupling between subsystems,  the following definition  formalizes the class of information decompositions considered in this paper. For each subsystem $i \in \nodes$, define the set
\begin{align*}
\Ncal(i) := \{j \in  \nodes_I^- (i) \ | \  \text{\eqref{eq:cond1}, \eqref{eq:cond2} are satisfied}   \},
\end{align*}
where the above conditions  are given by
\begin{align} 
\nodes_A^- (j) &\subseteq \nodes_I^- (i) , \label{eq:cond1} \\
\bigcup_{k \in \nodes_B^- (j)} \hspace{-.05in} \nodes_I^- (k) &\subseteq \nodes_I^- (i). \label{eq:cond2}
\end{align}
Condition \eqref{eq:cond1} requires that subsystem $i$ possess access to all states that directly affect subsystem $j$'s state through the system dynamics. Condition \eqref{eq:cond2} requires that subsystem $i$ have access to the local information of each subsystem whose control actions directly affect subsystem $j$'s state. This  ensures that subsystem $i$ is able to reconstruct all control actions that directly affect subsystem $j$'s state.
Collectively, conditions \eqref{eq:cond1} and \eqref{eq:cond2} can be interpreted as a requirement on the \emph{local nesting of information}, in the sense that if $j \in \Ncal(i)$, then  subsystem $i$ is assumed to  have access to all states and control actions that directly affect subsystem $j$'s state through the state equation. As a result, subsystem $i$ can explicitly reconstruct the local disturbance $w_j (t)$ acting on any subsystem $j \in \Ncal(i)$ based only on its local information $z_i (t+1)$ as follows:
\begin{align*}
\nonumber w_j (t) = x_j (t+1)  - \hspace{-.1in} \sum_{k \in \nodes_A^- (j)} A_{jk}(t) x_k(t)  - \hspace{-.1in}  \sum_{k \in \nodes_B^- (j)} B_{jk} (t) u_k (t) . 
\end{align*}

The local states of subsystems not belonging to $\Ncal(i)$ may contain information that can be influenced by prior control actions. We refer to these states as the \emph{information-coupling states} associated with subsystem $i$ at stage $t$, and  denote them by $x_{\couple(i)} (t)$ where  $$\couple (i) := \nodes_I^- (i) \setminus \Ncal(i).$$ 
The collection of information-coupling states across all subsystems are denoted by the  $x_\couple (t) \in \RR^{n_x^\couple}$, where $\couple := \bigcup_{i \in \nodes} \, \couple(i).$
The \emph{trajectory of information-coupling states} is denoted by
\begin{align*}
x_\couple := (x_\couple (0), \dots, x_\couple (T) ) \in \RR^{N_x^\couple},
\end{align*}
where $N_x^\couple := n_x^\couple (T+1).$ 
It will be notationally convenient to express the mapping from the state trajectory $x$ to its subvector $x_\couple$ in terms of the projection operator  $\proj: \RR^{N_x} \to \RR^{N_x^\couple}$, where $x_\couple = \Pi_\couple x$.

\begin{remar}[Partially Nested Information] It can be shown that the given information structure is \emph{partially nested} if and only if the set of  information coupling states is empty, i.e., $\couple = \emptyset$. 
\end{remar}

\subsection{Control Input Reparameterization}

The proposed information decomposition $ \nodes_I^- (i) = \Ncal(i) \cup \couple(i)$ suggests  a natural reparameterization of the control policy in terms of the following equivalent information set.

\begin{lemm}[Equivalent Information Sets] \label{lem:equi_info}
Define the information set $\zeta_i (t) := \{x_j^{0:t}  \, | \, j \in \couple (i) \} \cup \{w_j^{-1: t-1} \, | \, j \in \Ncal (i) \}$.
The sets $z_i (t)$ and $\zeta_i (t)$ are functions of each other for each subsystem $i$ and time $t$.
\end{lemm}

We omit the proof of Lemma \ref{lem:equi_info}, as it mirrors the proof of a closely related result in \cite[Lemma 1]{Lin2016}.
Lemma \ref{lem:equi_info} suggests the following equivalent parameterization of the local control input:
\begin{align}
u_i (t) = \surrpol_i (\pinfo_i (t), t), \label{eq:pol_equiv}
\end{align}
where $\surrpol_i (\cdot ,t)$ is a measurable function of its arguments. 
We let $\surrpol_i := ( \surrpol_i (\cdot, 0), .., \surrpol_i (\cdot, T-1) )$ and $\surrpol := (\surrpol_1, .., \surrpol_N)$ denote the reparameterized control policy for each subsystem $i \in \nodes$ and the full system, respectively.
With a slight abuse of notation, we express the input trajectory induced by the reparameterized  policy $\surrpol$ as $$u = \surrpol (w, x_\couple).$$
Finally, we denote by $\surrpolset$ the set of reparameterized decentralized control policies that respect the information constraints implied by Eq. \eqref{eq:pol_equiv}.
The reparameterization of the control policy  according to Eq. \eqref{eq:pol_equiv}  results in the following  equivalent reformulation of the original decentralized control problem \eqref{opt:decent}:
\begin{align}
\begin{alignedat}{8}
&\text{minimize} \quad && \EE \left[ x^\top R_x x + u ^\top R_u u \right] \\
& \text{subject to} \quad && \surrpol \in \surrpolset  \\
&&& \hspace{-.105in}\left. \begin{array}{l}
u = \surrpol (w, x_\couple) \\
x = B u + \wtoX w  \\
F_x x + F_u u + F_w w \leq g 
\end{array}
\hspace{-0.04in}\right\} \forall w \in \support.
\end{alignedat} \label{opt:decent_dFeedback}
\end{align}
Problem \eqref{opt:decent_dFeedback} is nonconvex, in general, if the set of information-coupling subsystems is nonempty, i.e.,  $\couple \neq \emptyset$. In Section \ref{sec:feas_control}, we construct a convex inner approximation of problem \eqref{opt:decent_dFeedback},  where the information-coupling states are assumed to behave as disturbances with bounded support, and the control policy is constrained in a manner that ensures consistency between  the assumed and actual behaviors of the information-coupling states.

\section{Decentralized Control Design via Contracts} \label{sec:feas_control}

In this section, we construct a convex inner approximation of the decentralized control design problem \eqref{opt:decent_dFeedback} 
by introducing a surrogate information structure in which the information-coupling states are modeled as \emph{fictitious disturbances}, which are  ``assumed" to take values in a given ``contract'' set.
To ``guarantee'' the satisfaction of this assumption, we impose a contractual constraint on the control policy requiring that the actual information-coupling states induced by the control policy belong to the contract set.
Given a fixed contract set, the resulting problem is a convex disturbance-feedback control design problem, whose feasible policies are shown  to be feasible for original problem \eqref{opt:decent_dFeedback}.

\subsection{Surrogate Information}
We associate a \emph{fictitious disturbance}  $\fic_i (t) \in \RR^{n_x^i}$ with each subsystem $i \in \nodes$ and time period $t = 0, \dots, T$. We let $\fic \in \RR^{N_x}$ denote the corresponding fictitious disturbance trajectory induced by these individual elements, which we model as a random vector \rtwo{whose distribution we take as given in this section.}\footnote{\rtwo{In Sec. \ref{sec:joint_opt}, we treat the probability distribution of the fictitious disturbance trajectory $\fic$ as a design choice. In Section \ref{sec:param_fict}, we propose a  parameterization of  the fictitious disturbance trajectory that facilitates the co-optimization of its distribution with the control policy via semidefinite programming.}}  The support of the fictitious disturbance trajectory, which we denote by $\contract \subset \RR^{N_x}$, is assumed to be a convex and compact set. We assume that the support of the pair of random vectors $(w, \fic)$ is given by  set $\support \times \contract$.

Treating the fictitious disturbances  as a surrogates for the information-coupling states, we define the \emph{surrogate local information} associated with subsystem $i$ as
\begin{align*}
\widetilde{\pinfo}_i (t) :=  \{v_j^{0:t} | j \in \couple (i) \} \cup \{w_j^{-1: t-1}  | j \in \Ncal (i) \}.
\end{align*}
Given a decentralized control policy $\surrpol \in \surrpolset$, the surrogate local information induces a \emph{surrogate control input} for each subsystem $i$ given by
$\widetilde{u}_i (t) := \surrpol_i ( \widetilde{\pinfo}_i (t), t).$
We denote the corresponding \emph{surrogate input trajectory}  
by
\begin{align*}
\widetilde{u} := \surrpol (w, \fic_\couple),
\end{align*}
where $\fic_\couple := \proj \fic$ is the projection of the fictitious disturbance $v$ onto the subspace associated with the information-coupling states.

\subsection{Surrogate Dynamics}

We replace the information-coupling states with  fictitious disturbances. This induces a \emph{surrogate system state} that evolves according to the following surrogate state equation:
\begin{align} 
\nonumber \widetilde{x}_i (t+1) = & \sum_{j \in \nodes \setminus \couple(i) } A_{ij} (t) \widetilde{x}_j (t) + \sum_{j \in \couple  (i)} A_{ij} (t) \fic_j (t)  \\
&  + \sum_{j=1}^N B_{ij} (t) \widetilde{u}_j (t) + w_i (t), \label{eq:surr_state_i}
\end{align}
where $\widetilde{x}_i (t)$ denotes the surrogate state of subsystem $i$ at time $t$. We require that the initial condition of the surrogate system equal that of the true system states agree, i.e.,  $\widetilde{x}_i (0) = x_i (0)$ for each subsystem $i$.
It will be convenient to express the  surrogate state dynamics in terms of the system trajectories according to
\begin{equation}
\widetilde{x} =  \Bsur \widetilde{u} +\Gsur w + \Hsur \fic_\couple. \label{eq:surrogate_state}
\end{equation}
\rtwo{It is straightforward to construct the lower block triangular matrices $\Bsur$, $\Gsur$, and $\Hsur$ given the problem data in \eqref{eq:surr_state_i}}

\subsection{Assume-Guarantee Contracts} \label{sec:assume-guarantee}

Under the proposed  surrogate information structure, the information-coupling states are assumed to behave as disturbances with bounded support.
In the following definition, we introduce a class of constraints on the family of admissible control policies, which guarantee that the information-coupling states behave in accordance with this assumption.

\begin{definitio}[Assume-Guarantee Contract] \label{def:assume-guarantee}
A control policy $\surrpol \in \surrpolset$ is said to satisfy the \emph{assume-guarantee contract} specified in terms of a \emph{contract set} $\contract_\couple \subseteq \RR^{N_x^\couple}$  if 
\ $\proj \widetilde{x} \in \contract_\couple$   for all  $(w , \fic_\couple  ) \in \Wcal \times \contract_\couple$, \
where $ \widetilde{x} =  \Bsur \surrpol (w,   \fic_\couple) + \Gsur w + \Hsur  \fic_\couple $.
\end{definitio}

Here, the set $\contract_\couple$ is referred to as a \emph{contract set}, as it specifies the set to which the information-coupling states are both assumed and required to belong.
The satisfaction of the assume-guarantee contract ensures that the \emph{surrogate} information-coupling states  $\widetilde{x}_\couple := \proj \widetilde{x}$ belong to the given contract set.
In the following lemma, we show that the \emph{actual} information-coupling states induced by the policy $u = \surrpol (w, x_\couple)$ are also guaranteed to belong to the contract set if the assume-guarantee contract is satisfied. 

\begin{lemm} \label{lem:assume-guarantee}
Let $\surrpol \in \surrpolset$ be a control policy that satisfies the assume-guarantee contract defined by the contract set $\contract_\couple \subseteq \RR^{N_x^\couple}$. It follows that $\Pi_\couple x\in \contract_\couple$ for all $w \in \support$, where $x = B \surrpol (w, x_\couple) + \wtoX w$.
\end{lemm}

Appendix \ref{pf:lem:assume-guarantee} contains the proof of Lemma \ref{lem:assume-guarantee}.
Prop. \ref{prop:feas} provides an inner approximation of the decentralized control design problem \eqref{opt:decent_dFeedback} by introducing an assume-guaranatee contractual constraint.

\begin{propositio} \label{prop:feas}
Let $\surrpol \in \Phi$ be a feasible control policy for the following problem:
\begin{align}
\begin{alignedat}{8}
&\text{\rm minimize} \quad && \EE \left[ \widetilde{x}^\top R_x \widetilde{x} + \widetilde{u} ^\top R_u \widetilde{u} \right] \\
& \text{\rm subject to} \quad && \surrpol \in \surrpolset \\
&&& \hspace{-.105in}\left. \begin{array}{l}
\widetilde{u} = \surrpol (w,  \fic_\couple) \\
\proj \widetilde{x} \in \contract_\couple \\
\widetilde{x} =  \Bsur \widetilde{u} +\Gsur w + \Hsur \fic_\couple \\
F_x \widetilde{x} + F_u \widetilde{u} + F_w w \leq g  
\end{array}
\hspace{-0.04in}\right\}  \forall (w, \fic_\couple) \in \support \times \contract_\couple.
\end{alignedat} \label{opt:contract}
\end{align}
It follows that $\phi$ is also  feasible for problem \eqref{opt:decent_dFeedback}.
\end{propositio}

We omit the proof of Prop. \ref{prop:feas}, as it is a direct consequence of Lemma \ref{lem:assume-guarantee}.
Given a fixed contract set $ \contract_\couple$, problem \eqref{opt:contract} is a convex disturbance feedback control design problem. The choice of the contract set  and fictitious disturbance distribution does, however, play an important role in determining the performance of the control policies that are optimal solutions to problem \eqref{opt:contract}. 
In Section \ref{sec:joint_opt}, we develop a systematic approach to enable the joint optimization of the contract set with the control policy using semidefinite programming.

\section{Policy-Contract Optimization} \label{sec:joint_opt}

In this section, we provide a semidefinite programming-based method to co-optimize the design of the decentralized control policy together with the contract set that constrains its design. 
As part of the proposed approach, we consider a restricted family of control policies that are affinely parameterized in both the disturbance and fictitious disturbance histories. We also parameterize the fictitious disturbance process as a causal affine function of a given (primitive) disturbance process---an approach that is similar in nature to the class of parameterizations that have been recently studied in the context of robust optimization with adjustable uncertainty sets \cite{zhang2017robust}. As one of our primary results in this section, we identify a structural condition on the family of allowable contract sets that facilitates the inner approximation of the resulting policy-contract optimization problem as a semidefinite program.

\subsection{Affine Decentralized Control Policies}
We consider affine disturbance-feedback  policies of the form
\begin{align}
\nonumber \widetilde{u}_i (t) =  u^o_i (t) &+ \sum_{ j \in \Ncal  (i) }   \sum_{s=-1}^{t-1} Q_{ij}^w (t,s+1) w_j (s)  \\
&+ \sum_{ j \in \couple (i)} \sum_{s=0}^t  Q_{ij}^\fic (t,s) \fic_j (s), \label{eq:aff1}
\end{align}
for $t = 0, .., T-1$ and $i = 1, .., N$. 
We refer to  $u^o_i (t) \in \RR^{n_u^i}$ as the open-loop control input, and the matrices $Q_{ij}^w (t,s+1) \in \RR^{n_u^i \times n_x^j}$ and $ Q_{ij}^\fic (t,s) \in \RR^{n_u^i \times n_x^j}$ as the feedback control gains.
The affine control policy \eqref{eq:aff1} can be expressed in terms of trajectories as
\begin{align}
\widetilde{u} = u^o +  Q^w  w + Q^\fic \fic, \label{eq:aff2}
\end{align}
where the  matrices $Q^w \in \RR^{N_u \times N_x}$ and $Q^\fic \in \RR^{N_u \times N_x}$ are both $T \times (T+ 1)$ block matrices, whose $(t, s)$-th blocks are defined as
\begin{align}
[Q^w (t,s) ]_{ij} &= \begin{cases}
Q_{ij}^w (t,s) & \text{if } j \in \Ncal (i) , \ t \geq s , \\
0 & \text{otherwise},
\end{cases}  \label{eq:Qw}\\
[Q^\fic (t,s)] _{ij} &= \begin{cases}
Q_{ij}^\fic (t,s) & \text{if } j \in \couple (i)  , \  t \geq s , \\
0 & \text{otherwise}.
\end{cases}  \label{eq:Qfic}
\end{align}
for $i = 1, \dots, N$ and  $j = 1, \dots, N$.
We let $\Qcal_N$ and $\Qcal_C$ denote the matrix subspaces  respecting the block sparsity patterns specified according to Eqs. \eqref{eq:Qw} and \eqref{eq:Qfic}, respectively.

\subsection{Affine Parameterization of the Fictitious Disturbance} \label{sec:param_fict}

We focus our analysis on fictitious disturbances that are expressed according to an affine transformation of a \emph{primitive disturbance}. Such a parameterization yields contract sets that have adjustable location, scale, and orientation.
Specifically, we let the random vector $\xi \in \RR^{N_x}$ denote the the \emph{primitive disturbance trajectory}, whose support $\Xi \subseteq \RR^{N_x}$ is assumed to be convex and compact. \rtwo{The distribution of the primitive disturbance trajectory $\xi$ is a design choice. In Assumption \ref{ass:ellip}, we propose a particular distribution that enables the co-optimization of the control policy and the contract set via semidefinite programming.}

We parameterize the fictitious disturbance trajectory affinely in the primitive disturbance as
\begin{align}
\fic := \overline{\fic} + \fictransform \xi. \label{eq:zeta}
\end{align}
The parameters  $\overline{\fic} \in \RR^{N_x}$ and $\fictransform \in \RR^{N_x \times N_x}$  control the shape of the resulting contract set, which is given by $\contract_\couple = \proj \left( \overline{\fic} \oplus \fictransform \Xi \right).$
Throughout the paper, we will restrict our attention to transformations \eqref{eq:zeta} in which the matrix parameter $\fictransform$ is both lower triangular and invertible. We denote the set of all such matrices by  $\Zcal \subset \RR^{N_x \times N_x}$.

The specification of the fictitious disturbance according to Eq. \eqref{eq:zeta} induces the following structure in the surrogate control input:
\begin{align}
\widetilde{u} = u^o + Q^\fic \overline{\fic} +  Q^w  w + Q^\fic \fictransform \xi. \label{eq:aff3}
\end{align}
We eliminate the bilinear terms in Eq. \eqref{eq:aff3} through the following  the change of variables:
\begin{align}
 \overline{u}:= u^o + Q^\fic \overline{\fic}   \ \text{ and } \ Q^\xi := Q^\fic \fictransform. \label{eq:change_of_var}
\end{align}
This change of variables gives rise to a reparameterization of the surrogate input trajectory as
\begin{align}
\widetilde{u} = \overline{u} +  Q^w  w + Q^\xi \xi, \label{eq:aff4}
\end{align}
where the matrix $Q^\xi \in \RR^{N_u \times N_x}$ must satisfy the sparsity constraint 
$Q^\xi \fictransform^{-1} \in \Qcal_C$
in order to ensure the satisfaction of the original sparsity constraint that $Q^\fic \in \Qcal_C$.
The affine parameterization of the control policy and the fictitious disturbance  according to   \eqref{eq:aff2} and \eqref{eq:zeta}, respectively, facilitates the co-optimization of the control policy and contract set as follows:
\begin{align}
\begin{alignedat}{8}
&\text{minimize} \quad && \EE \left[ \widetilde{x}^\top R_x \widetilde{x} + \widetilde{u} ^\top R_u \widetilde{u} \right] \\
& \text{subject to} \quad &&Q^w \in \Qcal_N, \,  Q^\xi\in \RR^{N_u \times N_x} , \,  \fictransform \in \Zcal \\
&&& \overline{u} \in \RR^{N_u}, \, \overline{\fic} \in \RR^{N_x} , \\
&&& Q^\xi \fictransform^{-1} \in \Qcal_C \\
&&& \hspace{-.105in} \left. \begin{array}{l}
\fic  =  \overline{\fic}  + \fictransform \xi\\
\widetilde{u} = \overline{u} +  Q^w  w + Q^\xi \xi \\
\widetilde{x} =  \Bsur \widetilde{u} +\Gsur w + \Hsur \fic_\couple  \\
\proj \widetilde{x} \in \proj  \left( \overline{\fic} \oplus \fictransform \support \right) \\
F_x \widetilde{x} + F_u \widetilde{u} + F_w w \leq g
\end{array}
\hspace{-0.04in}\right\} \forall (w, \xi) \in \support \times \Xi.
\end{alignedat} \label{opt:inner0_joint}
\end{align}

In contrast to problem \eqref{opt:contract}, where the contract set is considered fixed, problem \eqref{opt:inner0_joint}
allows for the co-optimization of the control  policy  together  with  the  contract set.
This increased flexibility in control design does,  however, result in the loss of problem convexity. Specifically,  the nonconvexity in  problem \eqref{opt:inner0_joint} stems from the assume-guarantee contractual constraint on the affine control policy, and the bilinear sparsity constraint  $Q^\xi \fictransform^{-1} \in \Qcal_C$. In Section \ref{sec:restrict}, we construct explicit convex inner approximations for each of these constraints, yielding a convex inner approximation of problem \eqref{opt:inner0_joint} in the form of a semidefinite program.

\subsection{Restricting the Contract Set}
\label{sec:restrict}

In what follows, we introduce additional restrictions on the set of allowable of matrix parameters $\mathcal{Z}$ that permit the reformulation of the nonconvex bilinear constraint $Q^\xi \fictransform^{-1} \in \Qcal_C$ as the linear constraint $Q^\xi \in \Qcal_C$.
\sam{We do so by imposing an additional subspace constraint on $\fictransform$ that guarantees the invariance of the subspace $\Qcal_C$ under multiplication with such matrices $\fictransform$.}
First, we impose a structural constraint the matrix parameter $\fictransform$ of the form
\begin{align}
\fictransform = \lambda I - Y. \label{eq:fictransform}
\end{align}
Here, $\lambda \geq 1$ is scalar parameter and $Y \in \RR^{N_x \times N_x}$ is a $(T+1) \times (T+1)$ strictly block lower triangular matrix of the form
\begin{align}
 Y = \bmat{0 \\ Y(1,0) & 0 \\ 
 \vdots &   \ddots  & \ddots\\
 Y (T, 0 ) &  \cdots & Y(T, T-1) & 0}. \label{eq:Y}
 \end{align}
In addition, each block of the matrix $Y$ is required to be a $N \times N$ block matrix, whose $(i, j)$-th block has dimension $n_x^i \times n_x^j$.   We impose an additional restriction on the structure of the matrix $Y$ in the form of sparsity constraints on each of its blocks that reflect the pattern of informational coupling between subsystems.
Specifically, we encode the pattern of informational coupling between subsystems according to a directed graph $\gcouple := (\nodes, \ecouple)$, where
\begin{align*}
\ecouple := \{(j, i) \in \edges_I \, | \, j \in \couple (i) \}.
\end{align*}
We let $\nodes_C^+ (i)$ denote the out-neighborhood of a node $i \in \nodes $ in the \emph{information coupling graph} $\gcouple$.
Using this graph, we impose a sparsity constraint  on each  block of the matrix $Y$ of the form:
\begin{align}
[Y (t, s) ]_{ij} = 0 \quad \text{if } \ \nodes_C^+ (i) \nsubseteq \nodes_C^+ (j)  \label{eq:subspace_Y}
\end{align}
for all $i,j \in \{1, \dots, N\}$, and $t, s \in  \{0,\dots, T\}$. We denote the subspace of all matrices $Y$ that respect these sparsity constraints by
\begin{align*}
    \Ycal (\gcouple) := \{ Y \in \RR^{N_x \times N_x} \ | \ Y \ \text{satisfies  \eqref{eq:Y} \& \eqref{eq:subspace_Y}} \}
\end{align*}

Lemma \ref{lem:QY} establishes the invariance of the subspace $\Qcal_C$ under multiplication by matrices 
$Y \in \Ycal (\gcouple)$. 
\rtwo{The result can be verified by direct calculation. The proof  is omitted due to space limitations.} 

\begin{lemm} \label{lem:QY} If $Q \in \Qcal_C$ and $Y \in \Ycal (\gcouple)$, then $Q Y \in \Qcal_C$.
\end{lemm}

The following useful result is a consequence of Lemma \ref{lem:QY}.

\begin{lemm} \label{lem:change_of_variable}
Let  $Y \in \Ycal (\gcouple)$  and $\lambda \in [1, \infty)$. It follows that $$\Qcal_C = \big\{ Q^\xi ( \lambda I - Y)^{-1} \, | \, Q^\xi \in \Qcal_C  \big\}.$$
\end{lemm}
\def \lmap{F}
\prf{Fix $\lambda \in [1, \infty)$ and $Y \in \Ycal (\gcouple)$.
Define the finite-dimensional linear map $\lmap: \RR^{N_u \times N_x} \to \RR^{N_u \times N_x}$ according to
$\lmap (Q) := Q (\lambda I - Y)^{-1}.$
We show that $\lmap (\Qcal_C  ) = \Qcal_C$. 
The invertibility of the linear map $\lmap$, in combination with the finite dimensionality of the subspace $\Qcal_C$, implies that $\lmap (\Qcal_C  ) = \Qcal_C$ if and only if $\lmap^{-1} (\Qcal_C ) \subseteq \Qcal_C$. Note that 
$\lmap^{-1} (Q) = Q (\lambda I - Y).$
The desired result follows, as Lemma \ref{lem:QY} implies that $QY \in \Qcal_C$ for each $Q \in \Qcal_C$ and $Y \in \Ycal (\gcouple)$.}

It follows from Lemma \ref{lem:change_of_variable} that if we impose the additional  constraint $$Z \in \{ \lambda I - Y \ | \ \lambda \in [1, \infty) \ \text{\&} \ Y \in \Ycal (\gcouple) \},$$ then the bilinear constraint $Q^\xi Z^{-1} \in \Qcal_C$ can be 
rewritten as an equivalent linear constraint $Q^\xi \in \Qcal_C$.

\subsection{Semidefinite Programming Approximation} \label{sec:sdp}

\rev{
In this section, we introduce a series of convex restrictions that culminate in the conservative approximation of
problem \eqref{opt:inner0_joint} as a semidefinite program. Our results rely on the following assumption concerning the distributions of the disturbance trajectory $w$ and primitive disturbance trajectory $\xi$.

\begin{assumptio} \label{ass:ellip}
Let $\Sigma \succ 0$. The disturbance trajectory $w$ is assumed to have an ellipsoidal support set given by
\begin{align*}
\support: = \big\{ z \in \RR^{N_x} \, \left|  \, z^\top \Sigma^{-1} z \leq 1\right. \big\}. 
\end{align*}
\rtwo{We also  assume that the primitive disturbance trajectory $\xi$ is an independent and identically distributed  copy of the disturbance trajectory $w$. In particular, this implies that $\Xi = \Wcal$.}
\end{assumptio}

In what follows, we employ  these assumptions to  reformulate the robust linear inequality constraints specified in problem \eqref{opt:inner0_joint} as second order cone constraints. We also construct a conservative approximation of the assume-guarantee contractual constraint in problem \eqref{opt:inner0_joint} according to a finite collection of linear matrix inequalities.
}

To lighten notation, we write the surrogate state trajectory $\widetilde{x}$ as
\begin{align*}
\widetilde{x} = \overline{x} +  \Xw  w +  \Xxi \xi ,
\end{align*}
where $\overline{x}  := \Bsur \overline{u} + \Hsur \proj \overline{\fic}$, $\Xw:= \Bsur Q^w + \Gsur$,  and $\Xxi:= \Bsur Q^\xi + \Hsur \proj (\lambda I - Y)$. 

The following result provides an equivalent reformulation of the robust linear inequality constraints specified in problem \eqref{opt:inner0_joint} as second-order cone constraints. The proof of Lemma \ref{lem:robust_linear} follows directly from the identity $\sup_{w \in \support} c^\top w = \lVert \Sigma^{1/2} c \rVert_2$ for all $c \in \RR^{N_x}$.
\begin{lemm} \label{lem:robust_linear} Let Assumption \ref{ass:ellip} hold.
The  robust  constraints 
$$F_x \widetilde{x} + F_u \widetilde{u} + F_w w \leq g \quad \forall  (w,  \xi) \in \support \times \Xi$$ are satisfied if and only if
\begin{align}
\nonumber & \left\lVert \Sigma^{1/2} e_i^\top  ( F_x \Xw +  F_u Q^w + F_w ) \right\rVert_2 + \left\lVert \Sigma^{1/2} e_i^\top ( F_x \Xxi + F_u Q^{\xi}) \right\rVert_2 \\
& \qquad \leq e_i^\top ( g - F_x \overline{x} - F_u \overline{u} )  \label{eq:robust_linear}
\end{align}
for $i = 1, \dots, m$, where $e_i \in \RR^m$ is the $i$th standard basis vector. 
\end{lemm}

We now address the nonconvexity that stems from the assume-guarantee contractual constraint specified in problem \eqref{opt:inner0_joint}: 
\begin{align} \label{eq:assgar}
    \proj \widetilde{x} \in \proj  \left( \overline{\fic} \oplus \fictransform \support \right) \quad \forall  (w,  \xi) \in \support \times \Xi.
\end{align}
First, notice  that, under the stated assumptions, the contractual constraint \eqref{eq:assgar} is equivalent to the following set containment constraint:
\begin{align}
\proj \left( \overline{x} \oplus \Xw \support \oplus  \Xxi \support \right)    \subseteq \proj  \left( \overline{\fic} \oplus  \fictransform \support  \right). \label{eq:containment}
\end{align}
The set containment constraint \eqref{eq:containment} requires that the Minkowski sum of two ellipsoids be contained within another ellipsoid.   We leverage on the following known result  to conservatively approximate this set containment constraint by a quadratic matrix inequality.

\begin{lemm}[Theorem 4.2 in \cite{durieu2001multi}] \label{lem:robust_quadratic}
Let Assumption \ref{ass:ellip} hold and let 
 $L_1, L_2, L_3 \in \RR^{m \times N_x}$, where it is assumed that $m \leq N_x$. If there exists a scalar $\alpha \in [0, 1]$ such that
\begin{align}
\alpha^{-1} L_1 \Sigma L_1^\top +  (1 - \alpha)^{-1} L_2 \Sigma L_2^\top \preceq  L_3 \Sigma L_3^\top. \label{eq:MI_lemma0}
\end{align}
then $L_1 \support \oplus L_2 \support \subseteq L_3 \support$.
\end{lemm}

\rev{
A direct application of Lemma \ref{lem:robust_quadratic} reveals that the set containment constraint \eqref{eq:containment} is satisfied if there exists a scalar $\alpha \in [0, 1]$ such that 
\begin{align}
&\proj (\overline{\fic} - \overline{x} ) = 0 \\
&\nonumber \proj \left( \alpha^{-1}  \Xw \Sigma \Xw^\top  + (1 - \alpha)^{-1} \Xxi \Sigma \Xxi^\top \right) \proj^\top  \\
& \qquad \preceq \proj (\lambda I - Y) \Sigma (\lambda I - Y)^\top \proj^\top. \label{eq:matrix_inequality1}
\end{align}
The quadratic matrix inequality \eqref{eq:matrix_inequality1} is  nonconvex in the decision variables $\lambda$ and $Y$. The following result provides a conservative approximation of the quadratic matrix inequality \eqref{eq:matrix_inequality1} as a linear matrix inequality.}  Appendix \ref{pf:lem:matrix_inequality} contains a proof of Lemma \ref{lem:matrix_inequality}.
\begin{lemm} \label{lem:matrix_inequality}
Let Assumption \ref{ass:ellip} hold. The set containment constraint \eqref{eq:containment} is satisfied if there exists a scalar $\beta \in [0, \lambda]$ such that
\begin{align}
&\proj \left(   \overline{x} - \overline{\fic} \right) =0 , \label{eq:LMI0}\\
&\bmat{\proj \widetilde{\Sigma} \proj^\top & \proj \Xw & \proj \Xxi \\
\Xw^\top \proj^\top & \beta \Sigma^{-1} & 0 \\
\Xxi ^\top \proj^\top & 0 & (\lambda - \beta) \Sigma^{-1} } \succeq 0, \label{eq:LMI}
\end{align}
where  $\widetilde{\Sigma} = \lambda \Sigma - Y \Sigma - \Sigma Y^\top.$
\end{lemm}

The combined application of Lemmas \ref{lem:change_of_variable}, \ref{lem:robust_linear}, and \ref{lem:matrix_inequality}
enables the conservative approximation of the nonconvex decentralized control design problem  \eqref{opt:inner0_joint} as a finite-dimensional semidefinite program. We state the resulting approximation in the following proposition.

\begin{propositio} \label{prop:final}
The  following semidefinite program is  a conservative approximation of problem \eqref{opt:inner0_joint}:
\begin{align}  \label{opt:inner4_joint}
\begin{alignedat}{8}
&\text{\rm minimize} \quad && {\rm Tr} \left( \Xxi^\top R_x  \Xxi M   +  \Xw^\top R_x  \Xw M \right) + \overline{x}^\top R_x  \overline{x}\\
&&& +  {\rm Tr} \left({Q^w}^\top R_u Q^w M  +   {Q^\xi}^\top R_u Q^\xi M  \right)+  \overline{u}^\top R_u  \overline{u} \\[3pt]
& \text{\rm subject to} \quad &&Q^w \in \Qcal_N,  \, Q^\xi \in \Qcal_C , \,  Y \in \Ycal (\gcouple) , \, \overline{u} \in \RR^{N_u}\\
&&& \Xw, \, \Xxi  \in \RR^{N_x \times N_x}, \ \overline{\fic}  , \,  \overline{x} \in \RR^{N_x} ,  \ \lambda , \,  \beta \in \RR_+ \\
&&& \lambda \geq \max \{ 1, \beta \} \\
&&& \overline{x} =  \Bsur \overline{u} + \Hsur \proj \overline{\fic} \\
&&& \Xw = \Bsur Q^w + \Gsur \\
&&& \Xxi = \Bsur Q^\xi + \Hsur \proj (\lambda I - Y) \\
&&& \eqref{eq:robust_linear}, \ \eqref{eq:LMI0}, \ \eqref{eq:LMI}.
\end{alignedat}  
\end{align}
The decision variables for problem \eqref{opt:inner4_joint} are the matrices $Q^w$, $Q^\xi$, $Y$, $\Xw$, $\Xxi$, the vectors $\overline{u}$, $\overline{v}$, $\overline{x}$, and the scalars $\lambda$ and $\beta$.\footnote{We note that the optimization variables $\overline{x}$, $\Xw$,  and $\Xxi$ have been introduced  to simplify the statement of problem \eqref{opt:inner4_joint}. These additional  variables can be eliminated via direct substitution. }
\end{propositio}

It follows from Prop. \ref{prop:final} that  any feasible solution to problem \eqref{opt:inner4_joint} can be mapped to an affine control policy that is guaranteed to be feasible for the original decentralized control design problem \eqref{opt:decent_dFeedback}, using the change of variables specified in  \eqref{eq:change_of_var}. Specifically, given any feasible solution to problem \eqref{opt:inner4_joint}, the following  affine control policy is guaranteed to be a feasible solution to problem \eqref{opt:decent_dFeedback}:
\begin{align}  \label{eq:final_aff}
u = \overline{u} \, + \, Q^w w \,+ \, Q^v(x- \overline{v}),
\end{align}
where $Q^v = Q^\xi (\lambda I - Y)^{-1}$. Here, $\overline{u} - Q^v \overline{v} \in \RR^{N_u}$ is the open-loop component of the control policy, and $Q^w \in \Qcal_N$ and  $Q^v \in \Qcal_C$  are the feedback control gains.  The satisfaction of the sparsity constraint $Q^v \in \Qcal_C$ is guaranteed by  Lemma \ref{lem:change_of_variable}.

\section{Illustrative Example} \label{sec:case_study}

\rtwo{In this section, we apply the control design methodology proposed in this paper to a linear time-invariant system composed of $N = 3$ subsystems, where the local state and input dimensions  of each subsystem are $n_x^i = n_u^i = 1$. The control horizon is set to $T = 15$.
The system matrices are given by
\begin{align*}
A(t) &= \bmat{0.5  & 0 & 0 \\
0.5 & 0.5 & 0 \\
0.5 & 1.2 & -1.2  }  \quad \text{and} \quad
B(t) = \diag (0.1, \,  1, \,  1) 
\end{align*}
for all $t = 0, \dots, T-1$. 
The system disturbance trajectory $w$ is assumed to be a uniformly distributed random vector over the hypersphere given by $\Wcal := \{w \in \RR^{3 (T+1)} \, \left| \, \lVert w \rVert_2 \leq 1  \right. \}$.
The system state and input trajectories are required to satisfy the constraints
\begin{align*}
\lVert x \rVert_{\infty} \leq x_{\max} \quad \text{and} \quad \lVert u \rVert_{\infty} \leq 2.5 \quad  \forall w \in \Wcal,
\end{align*}
where  $x_{\max} \in \RR_+$ is constraint parameter that we will vary in our numerical studies. The cost matrices are specified as
\begin{align*}
R_x = I_{T+1} \otimes \diag (0.1, \, 0.1, \, 2) \quad \text{and} \quad
R_u = I_{T} \otimes \diag (5,\, 5, \, 1),
\end{align*}
where $\otimes$ denotes the Kronecker product operator.}

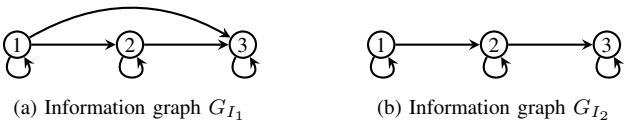
\begin{figure}[http]
\tikzstyle{regular_node} = [circle, draw = black, thick,  fill = white,  inner sep = 0pt, minimum size = 3.5mm]

\tikzstyle{edge_dir} = [->,>=stealth, thick] 
\tikzstyle{edge_bidir} = [latex-latex]
\tikzstyle{loop_spec} = [out = 240, in = -60, looseness = 5.5]

\centering

\begin{minipage}[b]{0.45\linewidth}
\centering
\begin{tikzpicture}

	\node at (-3,0) [regular_node] (n1)   {\footnotesize 1};
	\node at (-1.5,0) [regular_node] (n2) {\footnotesize 2};
	\node at (0,0) [regular_node] (n3)    {\footnotesize 3};
	
	\draw [edge_dir] (n1) to (n2);
	\draw [edge_dir] (n2) to (n3);

	\draw [edge_dir] (n1) to [bend left = 30] (n3);
	\draw [edge_dir] (n1) to [loop_spec] (n1);
 	\draw [edge_dir] (n2) to [loop_spec] (n2);
 	\draw [edge_dir] (n3) to [loop_spec] (n3);

\end{tikzpicture}

{\footnotesize  (a) Information graph $\gr_{I_1}$}
\end{minipage}
\hspace{2em}
\begin{minipage}[b]{0.45\linewidth}
\centering
\begin{tikzpicture}

	\node at (-3,0) [regular_node] (n1)  {\footnotesize 1};
	\node at (-1.5,0) [regular_node] (n2) {\footnotesize 2};
	\node at (0,0) [regular_node] (n3)  {\footnotesize 3};

	\draw [edge_dir] (n1) to (n2);
	\draw [edge_dir] (n2) to (n3);
 	\draw [edge_dir] (n1) to [loop_spec] (n1);
 	\draw [edge_dir] (n2) to [loop_spec] (n2);
 	\draw [edge_dir] (n3) to [loop_spec] (n3);

\end{tikzpicture}

{\footnotesize (b) Information graph $\gr_{I_2}$}
\end{minipage}

\caption{\footnotesize Two  different information graphs examined in this case study.}

\label{fig:info_graph}

\end{figure}

\rtwo{In this case study, we evaluate the performance our proposed control design method using two different information graphs depicted in Fig. \ref{fig:info_graph}. The information graph $\gr_{I_1}$ in Fig. \ref{fig:info_graph}(a) induces a  decentralized control design problem with a \emph{partially nested} information structure, while  the information graph $\gr_{I_2}$ in Fig. \ref{fig:info_graph}(b) results in a \emph{nonclassical} information structure. We  also compare our approach to a  closely related contract-based decentralized control design method  from  \cite{darivianakis2018decentralized}.}

\vspace{-.1in}

\subsection{Numerical Results and Discussion}

\begin{figure}[h]
    \centering
    \includegraphics[width = 0.95 \linewidth]{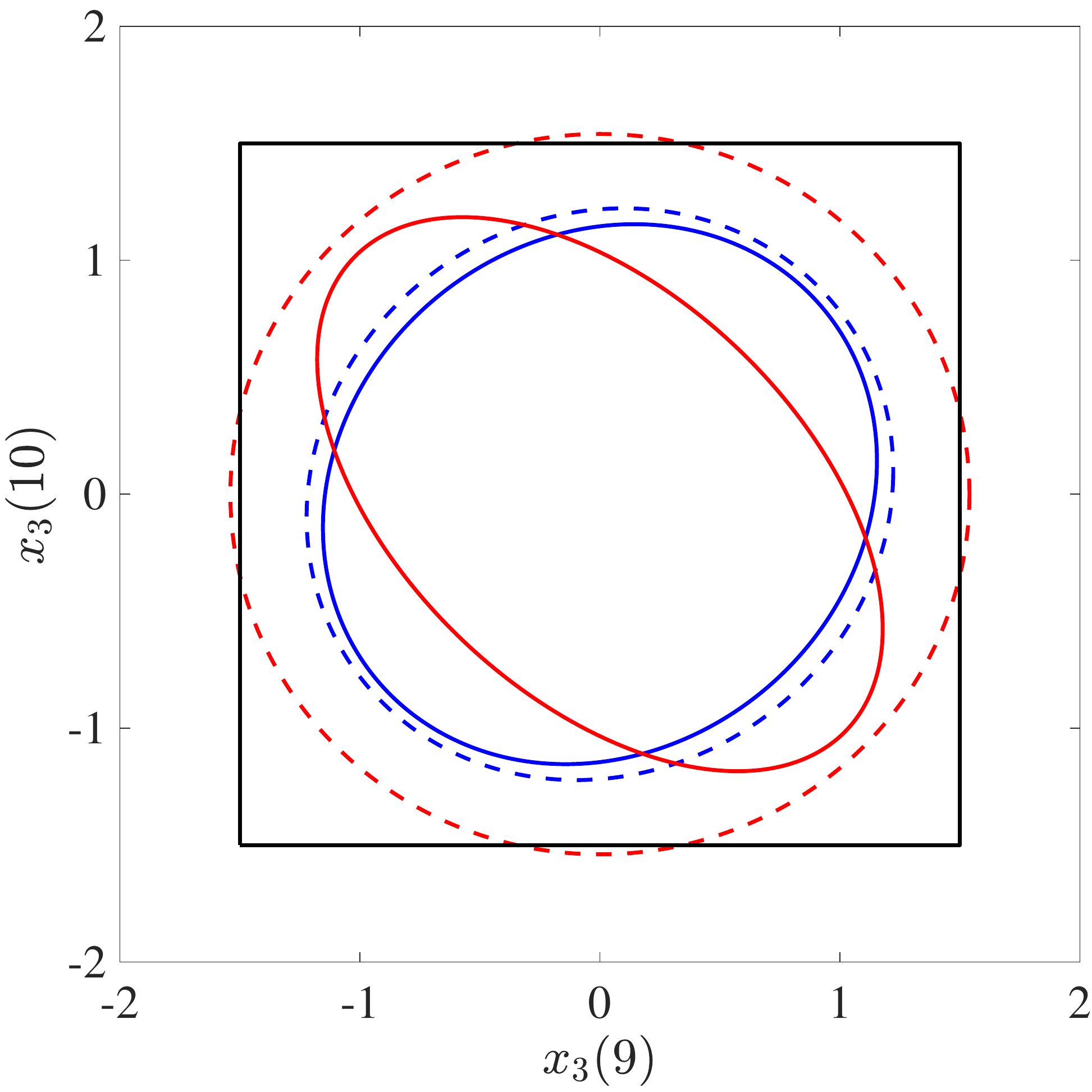}
    \caption{\footnotesize \rtwo{Plots of the state-support sets (solid boundaries) and the contract sets (dashed boundaries) associated with states $(x_3 (9), x_3 (10))$. The sets  generated by the method in \cite{darivianakis2018decentralized} using the  information graph $\gr_{I_1}$ are depicted in red. The sets  generated by our method using the  information graph $\gr_{I_2}$ are depicted in blue. The state constraint set is depicted in black.}}
    \label{fig:support}
\end{figure}

\rtwo{Focusing initially on the class of decentralized control problems induced by the information graph  $\gr_{I_1}$ in Fig. \ref{fig:info_graph}(a), we compare the performance of our control design approach against that of the method proposed in \cite{darivianakis2018decentralized}. In Fig. \ref{fig:parametric}, we plot the cost incurred by controllers computed based on our method (solid blue line)  and  the cost incurred by controllers computed according to the method in  \cite{darivianakis2018decentralized} (solid red line) 
for different values of the state constraint parameter $x_{\max}$ ranging between 1 and 2.5.  We also plot a lower bound (dashed black line) on the optimal value of each problem instance \eqref{opt:decent} using a recently proposed convex relaxation \cite{lin2019convexinfo}. 
Notice that our control design technique results in  controllers that  strictly outperform the controllers generated by the method from \cite{darivianakis2018decentralized} for each value of the constraint parameter $x_{\max}$.    Moreover, our method results in controllers that appear to be globally optimal for constraint parameter values $x_{\max} \in  [1.3, \, 2.5]$, as they achieve the lower bound on the optimal cost in that parameter regime.  The cost curves are not plotted for values of the constraint parameter $x_{\max}$ where the underlying control design method fails to return a feasible controller.  Notice that the method from this paper is more successful in generating  feasible controllers for values of the constraint parameter $x_{\max}$ close to one.

\begin{figure}[http]
\centering
\includegraphics[width = 1 \linewidth]{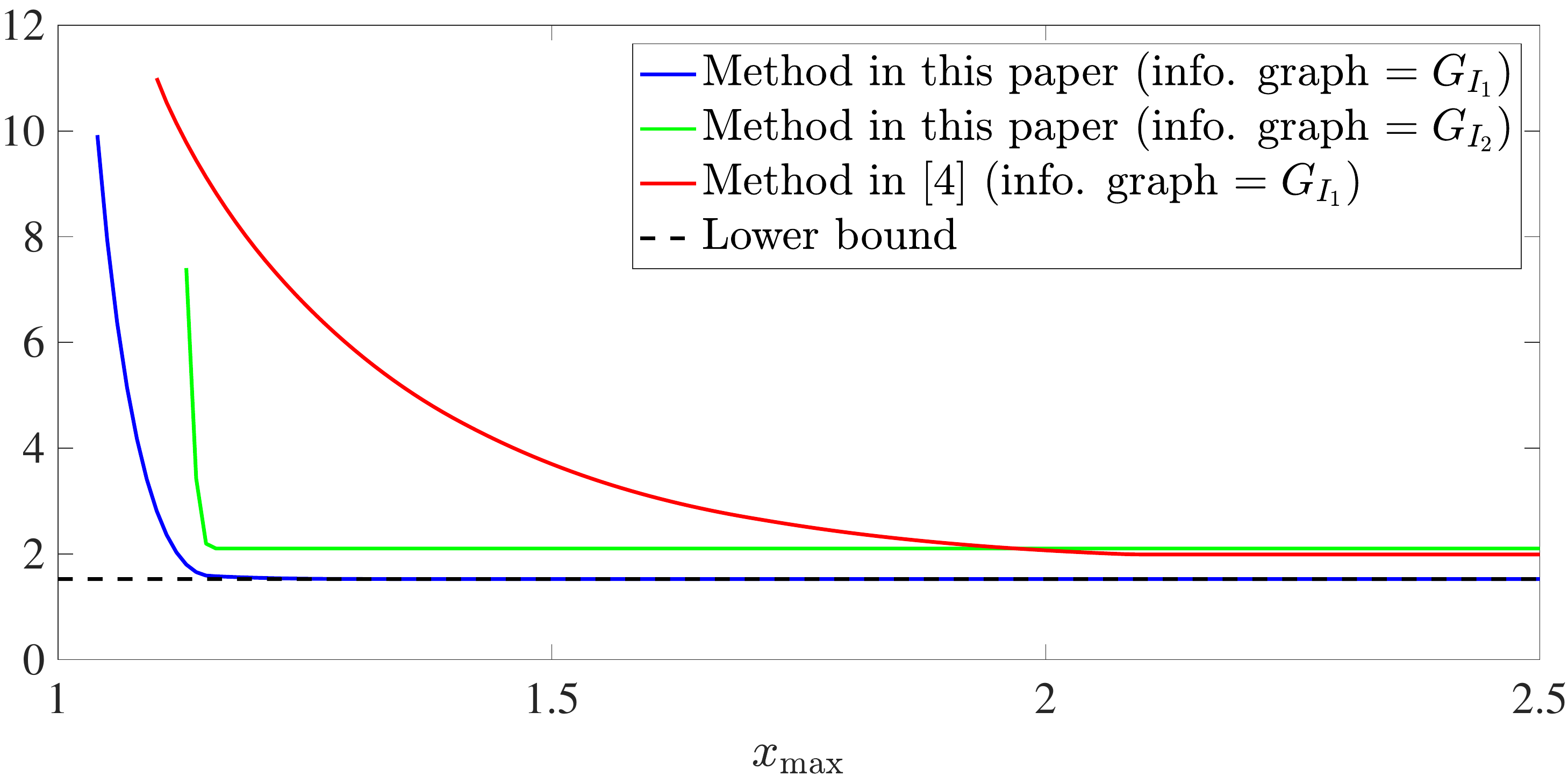}
\caption{\footnotesize \rtwo{For different values of the constraint parameter $x_{\max}$ ranging between 1 and 2.5, we plot the cost incurred by controllers generated by: (i) our method using the information graph $\gr_{I_1}$ (blue solid line); (ii) our method using the information graph $\gr_{I_2}$ (green solid line); (iii) the method in \cite{darivianakis2018decentralized} using the information graph $\gr_{I_1}$ (red solid line). The cost curves are not plotted for values of the constraint parameter $x_{\max}$ where the underlying control design method fails to return a feasible controller. The black dashed line depicts a lower bound on the optimal value of each problem instance using the convex relaxation method from \cite{lin2019convexinfo}.}}

\label{fig:parametric}
\end{figure}

In Fig. \ref{fig:parametric}, we also  plot the cost incurred by controllers generated by our method (solid green line) using the smaller information graph $\gr_{I_2}$ depicted in Fig. \ref{fig:info_graph}(b).
Despite using less information, our approach still results in controllers that substantially outperform those generated by \cite{darivianakis2018decentralized} for a wide range of constraint parameter values. 

There are two possible explanations for these observed differences in performance between the two methods. First, the approximation technique proposed in \cite{darivianakis2018decentralized}  treats all neighboring states to a subsystem as fictitious disturbances, while our approach only treats neighboring states that result in  `information coupling' as fictitious disturbances. Second, the method from \cite{darivianakis2018decentralized} does not permit the rotation of the primitive contract set when co-optimizing its specification with the decentralized control policy.  Because of this limitation, the method in \cite{darivianakis2018decentralized} may fail to generate contract sets that accurately capture the spatial and inter-temporal correlation between different states, resulting in more conservative approximations. 
We illustrate this limitation in Fig. \ref{fig:support}, where we plot the state-support sets and contract sets generated by each method for $x_3(t)$ at two consecutive time periods $t=9$ and $t=10$. Notice that our method generates an ellipsoidal contract set that more accurately captures 
the correlation between $x_3(9)$ and $x_3(10)$.
}

\appendix

\subsection{Proof of Lemma \ref{lem:assume-guarantee}} \label{pf:lem:assume-guarantee}

\def\projt{\Pi^{0:t}}

We require several  definitions. Define the projection operators 
\begin{align*}
\projt_\couple &:= \bmat{I_{n_x^\couple (t+1)} & 0_{n_x^\couple (t+1) \times n_x^\couple (T - t)}} \proj
\end{align*}
for  $t = 0, .., T$.
It follows $x_\couple^{0:t} = \projt_\couple x$ for each time $t \in \{0, .., T\}$. Additionally, we define the function $f_\surrpol: \RR^{N_x} \times \RR^{N_x^\couple} \to \RR^{N_x}$ as
\begin{align}
f_\surrpol (w, \fic_\couple ) := \Bsur \surrpol (w,   \fic_\couple) + \Gsur w + \Hsur  \fic_\couple. \label{eq:f_phi}
\end{align}
\rev{The function $f_\surrpol$ maps  $(w,\fic_\couple)$ to the surrogate system state $\widetilde{x}$  under the control policy $\surrpol$.}
Lemma \ref{lem:f_phi} shows that the function $f_\surrpol$ is \emph{strictly causal} in $\fic_\couple$, and that $x = f_\surrpol (w, x_\couple)$.
We state it without proof, as it  follows from the  structure of the  surrogate state equation \eqref{eq:surr_state_i}.

\begin{lemm} \label{lem:f_phi}
Let $\surrpol \in \surrpolset$ and $x = B\surrpol (w, x_\couple) + \wtoX w$. The following properties hold.
\begin{enumerate}[(i)]\setlength{\itemsep}{.05in}
\item Let $t \in \{0, \dots, T\}$ and  $w \in \RR^{N_x}$. It holds that 
$$\Pi^{0:t-1}_\couple (v - v') = 0 \ \Longrightarrow \ \projt_\couple f_\surrpol (w, v_\couple) = \projt_\couple f_\surrpol (w, v_\couple').$$
\item The state trajectory satisfies  $x = f_\surrpol (w, x_\couple)$ for all $w \in \RR^{N_x}$.
\end{enumerate}
\end{lemm} 

\vspace{.03in}

Now, fix $w \in \support$ and let $x = B\surrpol (w, x_\couple) + \wtoX w$. \rev{To complete the proof, it suffices to show that $x_\couple^{0:t} \in \projt_\couple \proj^\top \contract_\couple$ for all $t \in \{0, \dots, T\}.$  In particular, the satisfaction of this condition for $t = T$ implies that $x_\couple \in \proj \proj^\top \contract_\couple = \contract_\couple$. We prove this by induction in $t$.}

\emph{(Base step).} \, We have that $x_\couple (0) = \widetilde{x}_\couple (0) \in \Pi_\couple^{0:0} \proj^\top \contract_\couple$, where the first equality follows from the initial condition of the surrogate state equation \eqref{eq:surr_state_i}.

\emph{(Induction step).} \, Assume that $x^{0:t-1}_\couple \in \Pi^{0:t-1}_\couple \proj^\top \contract_\couple$. We complete the proof by showing that $x^{0:t}_\couple \in \projt_\couple \proj^\top \contract_\couple$. Fix a $v \in \RR^{N_x}$ that satisfies $\proj v \in \contract_\couple$ and $x^{0:t-1}_\couple = \Pi_\couple^{0:t-1} v$ (which is guaranteed to exist by our induction hypothesis). 
We have that
\begin{align*}
x_\couple^{0:t} &=  \projt_\couple   x =  \projt_\couple   f_\surrpol (w, x_\couple) = \projt_\couple   f_\surrpol (w, v_\couple)\\
&= \projt_\couple \proj^\top \proj  f_\surrpol (w, v_\couple) \in \projt_\couple \proj^\top \contract_\couple.
\end{align*}
The second equality follows from property (ii) in Lemma \ref{lem:f_phi}; the third equality follows from a combination of property (i) in Lemma \ref{lem:f_phi} and the equality $x^{0:t-1}_\couple = \Pi_\couple^{0:t-1} v$; the fourth equality follows from the identity that $\projt_\couple \proj^\top \proj = \projt_\couple$ for each time $t$; and the final inclusion  follows from the assumption that the policy $\surrpol$ satisfies the assume-guarantee contract specified by the contract set $\contract_\couple$.

\subsection{Proof of Lemma \ref{lem:matrix_inequality}} \label{pf:lem:matrix_inequality}

It suffices to show that the matrix inequality \eqref{eq:matrix_inequality1} is satisfied if the LMI \eqref{eq:LMI} is satisfied.
Define $\beta := \alpha \lambda$, and divide both sides of the matrix inequality \eqref{eq:matrix_inequality1} by $\lambda$. It follows that the matrix inequality \eqref{eq:matrix_inequality1} is satisfied if and only if there exists $\beta \in [0, \lambda]$, such that
\begin{align}
&\nonumber \beta^{-1} \proj \Xw \Sigma \Xw^\top \proj^\top + (\lambda - \beta)^{-1} \proj \Xxi \Sigma \Xxi^\top \proj^\top  \\
& \qquad \preceq \proj (\lambda \Sigma - Y \Sigma - \Sigma Y^\top + \lambda^{-1} Y \Sigma Y^\top) \proj^\top. \label{eq:matrix_inequality2}
\end{align}
Since $\lambda^{-1} Y \Sigma Y^\top \succeq 0$, the matrix inequality \eqref{eq:matrix_inequality2} is satisfied if 
\begin{align*}
    \beta^{-1} \proj \Xw \Sigma \Xw^\top \proj^\top + (\lambda - \beta)^{-1} \proj \Xxi \Sigma \Xxi^\top \proj^\top  \preceq \proj \widetilde{\Sigma} \proj^\top,
\end{align*}
where $\widetilde{\Sigma} := \lambda \Sigma - Y \Sigma - \Sigma Y^\top$. It follows from Schur's Lemma that the above matrix inequality is satisfied if and only if the LMI \eqref{eq:LMI} is satisfied. This completes the proof.


\bibliographystyle{abbrv}
\bibliography{decentralized_bib}{\markboth{References}{References}}

\end{document}